\documentclass[11pt,reqno]{amsart}

\usepackage[T1]{fontenc}

\usepackage{amsmath,amsfonts,amssymb,amsthm}
\usepackage{enumerate}
\usepackage{graphicx}

\usepackage{url} 
\usepackage[colorlinks=true, linkcolor=blue, citecolor=blue, urlcolor=black]{hyperref}

\numberwithin{equation}{section}
\newtheorem{theorem}{Theorem}

\newtheorem{lemma}{Lemma}
\newtheorem{definition}{Definition}
\newtheorem{corollary}{Corollary}

\theoremstyle{remark}
\newtheorem*{remark}{Remark}

\newcommand{\essinf}{\operatornamewithlimits{ess\,inf}}
\newcommand{\esssup}{\operatornamewithlimits{ess\,sup}}

\allowdisplaybreaks

\begin{document}

\title[Approximation in variable Lebesgue spaces]
{Approximation by classical summability methods in variable Lebesgue spaces}

\author{Tengiz Kopaliani}
\address{Tengiz Kopaliani \\
    Department of Mathematics\\
	Faculty of Exact and Natural Sciences\\
	Tbilisi State University\\
	1, Ilia Chavchavadze Av., Tbilisi, 0179, Georgia}
\email{tengiz.kopaliani@tsu.ge}

\author{Alina Shalukhina}
\address{Alina Shalukhina \\
	Center for Mathematics and Applications and Department of Mathematics \\
	NOVA School of Science and Technology, NOVA University Lisbon\\
	Quinta da Torre, Caparica, 2829--516, Portugal}
\email{a.shalukhina@campus.fct.unl.pt}

\author{George Tephnadze}
\address{George Tephnadze \\
	School of Science and Technology \\
	The University of Georgia \\
	77a, Merab Kostava St., Tbilisi, 0171, Georgia}
\email{g.tephnadze@ug.edu.ge}

\thanks{The second author is funded by national funds 
through the FCT -- Funda\c{c}\~{a}o para a Ciência e a Tecnologia, I.P., 
under the scope of the projects UID/00297/2025 and UID/PRR/00297/2025 
(Center for Mathematics and Applications – NOVA Math) and under the scope of 
the PhD scholarship UI/BD/154284/2022.
The third author is supported by Shota Rustaveli National Science Foundation 
grant no. FR-24-698.}



\begin{abstract}
We investigate approximation by Fej\'er and N\"orlund means
with respect to the Walsh system in variable Lebesgue spaces. 
We prove that the Fej\'er means are uniformly 
bounded on the weighted space $L^1_w[0,1)$ if and only if $w\in A_1(\mathcal{D})$. 
We then apply Rubio de Francia extrapolation to show that both means are uniformly 
bounded on $L^{p(\cdot)}[0,1)$ with $1\le p_-\le p_+<\infty$ if and only if 
$p(\cdot)\in\mathcal{A}(\mathcal{D})$. Finally, we estimate the 
approximation errors via the modulus of smoothness and deduce 
convergence rates for variable Lipschitz classes.
\end{abstract}

\maketitle

\section{Introduction}

Fej\'er means and, more generally, N\"orlund means provide 
well-established techniques for function approximation. For the classical 
trigonometric system, the Fej\'er means $\sigma_nf$ of a $2\pi$-periodic function
$f\in L^1[-\pi,\pi)$ are defined via the convolution
\[
\sigma_n f(x):=(f* K_n)(x)
=\frac1{2\pi}\int_{-\pi}^\pi f(t)\,K_n(x-t)\,dt,
\]
where the trigonometric Fej\'er kernel
\[
K_n(x):=\frac1n\,
\frac{\sin^2(nx/2)}{\sin^2(x/2)}
\] 
is bounded by a radially decreasing, integrable majorant. 
This bound implies that the Fej\'er means are uniformly, pointwise dominated by the Hardy--Littlewood maximal function~$Mf$: one has the estimate
\begin{equation}\label{eq:sigma-M}
|\sigma_nf(x)|\le C\,Mf(x)
\end{equation}
with a constant $C>0$ independent of $n$ and $f$; see Hardy and 
Littlewood~\cite[Theorem~20]{HL30} and, e.g., Grafakos~\cite[Theorem~2.1.10]{cLG14}.

Estimate~\eqref{eq:sigma-M} ensures that the boundedness
of the maximal operator~$M$ on a given normed lattice implies the uniform boundedness
of~$\sigma_n$. A standard density argument then yields the norm 
convergence $\sigma_nf\to f$. For example, the operators $\sigma_n$ are uniformly 
bounded on the weighted Lebesgue spaces $L^p_w[-\pi,\pi)$ for $1<p<\infty$ if 
the weight $w$ belongs to the Muckenhoupt $A_p$ class~\cite[Theorem~7.1.9]{cLG14}. 
Furthermore, while $M$ fails to be bounded on $L^1_w[-\pi,\pi)$, the uniform 
boundedness of $\sigma_n$ still holds when $w \in A_1$. This is established by passing 
the convolution to the weight via Fubini's theorem: combining the pointwise bound 
$\sigma_nw\le C\, Mw$ with the $A_1$ condition $Mw \le C\,w$ results in 
$\|\sigma_nf\|_{L^1_w}\le C\,\|f\|_{L^1_w}$.

For $1<p<\infty$, an exact equivalence holds for the Fej\'er means with respect to 
the \emph{Walsh system}. As a consequence of Kazarian's result~\cite[Theorem~4]{KK82}, 
the Walsh--Fej\'er operators $\sigma_n$ are uniformly bounded on $L^p_w[0,1)$ if 
and only if $w$ belongs to the dyadic Muckenhoupt class $A_p(\mathcal{D})$. However, 
the boundary case $p=1$ has remained open. Resolving this gap is the first 
contribution of our paper (Theorem~\ref{th:Fejer-A1}): 
\begin{center}
\emph{The Walsh--Fej\'er 
operators~$\sigma_n$ are uniformly bounded on $L^1_w[0,1)$\\ if and only if 
$w \in A_1(\mathcal{D})$.} 
\end{center}

Because the Walsh--Fej\'er kernel $K_n$ lacks a decreasing, integrable majorant, the 
simple argument relying on~\eqref{eq:sigma-M} is invalid. Our proof is therefore 
inherently dyadic and rests on bitwise shifts of dyadic intervals.

Building on this endpoint result, we apply a version of the Rubio de Francia 
extrapolation theorem for variable Lebesgue spaces to establish our second 
contribution (Theorem~\ref{th:main-variable}):
\begin{center}
\emph{
Given a variable exponent $p(\cdot)$ satisfying $1\le p_-\le p_+<\infty$,\\
the Walsh--Fej\'er operators $\sigma_n$ are uniformly bounded 
on $L^{p(\cdot)}[0,1)$\\
if and only if $p(\cdot)$ is in the dyadic 
Muckenhoupt-type class $\mathcal{A}(\mathcal{D})$, i.e.,
\[
\sup_{I\in\mathcal{D}}
\frac{1}{|I|}\,\|\chi_I\|_{p(\cdot)} \|\chi_I\|_{p'(\cdot)}<\infty
\]
where the supremum is taken over all dyadic intervals~$I$ on $[0,1)$. 
}
\end{center}

Our theorem sharpens the earlier result by 
Weisz~\cite[Corollary~8.7]{FW22}, who established the norm convergence 
$\|\sigma_nf-f\|_{p(\cdot)}\to0$ under the restrictive assumption 
$1<p_-\le p_+<\infty$ and the standard log-H\"older continuity condition
\[
|p(x)-p(y)|\le\frac{C}{\log (e+1/\vert x-y\vert)}
\quad\textnormal{for all }x,y\in[0,1).
\]
Theorem~\ref{th:main-variable} advances the 
current literature in three ways:
\begin{enumerate}
\item[\textnormal{(a)}] \emph{Inclusion of the endpoint:} We cover the boundary case 
$p_-=1$, which is frequently excluded in variable Lebesgue space theory.

\item[\textnormal{(b)}] \emph{Sharp assumptions:} We replace 
log-H\"older continuity with the condition~$p(\cdot)\in\mathcal{A}(\mathcal{D})$. 
While the former is merely sufficient for the boundedness of the maximal 
operator~\cite[Theorem~4.3.8]{DHHR11}, $\mathcal{A}(\mathcal{D})$ fully 
characterizes the boundedness of the dyadic maximal operator $M^\mathcal{D}$ 
on $L^{p(\cdot)}[0,1)$, even in the extended regime $1<p_-\le p_+\le\infty$ 
(Theorem~\ref{th:Kopaliani}).

\item[\textnormal{(c)}] \emph{Exact equivalence for the Walsh system:} 
As a result, $p(\cdot)\in\mathcal{A}(\mathcal{D})$ is 
necessary and sufficient for the uniform boundedness of the 
Walsh--Fej\'er operators~$\sigma_n$. Sufficiency follows from the maximal operator 
bounds, whereas necessity relies on the dyadic 
structure of the Walsh system.
\end{enumerate}
Note that extending this characterization to the case $p_+=\infty$ remains an open 
problem. Since Rubio de Francia extrapolation does not apply here, 
analyzing the Walsh--Fej\'er means requires alternative methods.

Not only do we establish the uniform boundedness of $\sigma_n$ and its 
equivalence to the norm convergence $\|\sigma_nf-f\|_{p(\cdot)}\to0$, but we also 
estimate the approximation error via the dyadic modulus of smoothness
\begin{equation}\label{eq:mod-smooth}
\omega_{p(\cdot)}(f,2^{-k}):=\|f-S_{2^k}f\|_{p(\cdot)},
\end{equation}
where $S_{2^k}f$ is the $2^k$-th dyadic partial sum of the Walsh--Fourier 
series of $f$. Because variable Lebesgue spaces lack translation invariance, 
the modulus~\eqref{eq:mod-smooth} naturally replaces the classical 
modulus of continuity. While this modulus of smoothness is widely used for 
trigonometric approximation in 
$L^{p(\cdot)}$ spaces, we apply it here for the first time to bound the approximation 
error with respect to the Walsh system. Establishing this error estimate is 
our third contribution (Theorem~\ref{th:Fej-var-Lp}):
\begin{center}
\emph{If $p(\cdot)\in\mathcal{A}(\mathcal{D})$ satisfies $1\le p_-\le p_+<\infty$,\\ 
there exists $C>0$ such that for all $f\in L^{p(\cdot)}[0,1)$ and 
natural numbers $n$,
\[
\|\sigma_n f-f\|_{p(\cdot)}
\le C\sum_{k=0}^{|n|} 2^{k-|n|}\,\omega_{p(\cdot)} (f,2^{-k}),
\]
where the integer $|n|$ is defined by $2^{|n|}\le n<2^{|n|+1}$.
} 
\end{center}

Furthermore, this estimate yields rates of convergence 
$\|\sigma_nf-f\|_{p(\cdot)}\to0$ for functions $f$ belonging to the variable 
Lipschitz classes $\textnormal{Lip}(\alpha,p(\cdot))$,
as detailed in Corollary~\ref{cor:Lip-Fejer}. For the trigonometric system,
similar estimates were obtained by Akg\"un~\cite{RA11}. 

Finally, we apply our three main results for the Walsh--Fej\'er means to
derive their exact analogues for regular N\"orlund means~$t_n$ in 
Theorems~\ref{th:Norlund-Ap}, \ref{th:Norlund-var}, and~\ref{th:Norlund-est}.
This extends the investigation of N\"orlund and Riesz means in variable Lebesgue 
spaces, which was conducted in the trigonometric setting by Guven and 
Israfilov~\cite{GI10}. In a broader context, our findings build upon and extend the 
modern framework of dyadic harmonic analysis developed in the recent monograph by 
Persson, Tephnadze, and Weisz~\cite{PTW22}.

\section{Basics of Walsh--Fourier Analysis}

Let $\mathbb{P}$ denote the set of positive integers and 
$\mathbb{N}:=\mathbb{P}\cup\{0\}$. 
Every $n\in\mathbb{N}$ and $x\in[0,1)$ can be uniquely expressed via the 
dyadic expansions 
\[
n=\sum_{i=0}^{\infty}n_i 2^i
\quad\text{and}\quad
x=\sum_{i=0}^{\infty}\frac{x_i}{2^{i+1}},
\]
where $n_i,x_i\in\{0,1\}$, with $x_i$ referred to as the $i$-th dyadic coordinate 
of~$x$. To ensure uniqueness of the representation for $x$, we 
exclude expansions where $x_i=1$ for all sufficiently large $i$.
Define the \emph{dyadic addition} $\dotplus$ on the interval $[0,1)$ by
\[
x\dotplus y=\sum_{i=0}^{\infty}\frac{(x_i+y_i) \bmod 2}{2^{i+1}}.
\]

By a \emph{dyadic interval}, we mean an interval of the form 
$[k\,2^{-n},(k+1)\,2^{-n})$ for some $n,k\in\mathbb{N}$ with $0\le k<2^n$. We denote 
the collection of all such intervals by $\mathcal{D}$. For any $x\in[0,1)$ and 
$n\in\mathbb{N}$, the unique dyadic interval of length $2^{-n}$ containing $x$ 
is given by
\[
I_n(x):=\{y\in[0,1)\,:\,y_0=x_0,\ldots,y_{n-1}=x_{n-1}\},
\]
where $I_0(x):=[0,1)$. For brevity, we write $I_n:=I_n(0)=[0,2^{-n})$.

Let $\{w_n\}_{n\in \mathbb{N}}$ represent the \emph{Walsh system} in Paley ordering, 
defined for $x\in[0,1)$ by 
\[
w_n(x):=\prod_{i=0}^{\infty}\,(-1)^{n_i x_i}.
\]
For $f\in L^1[0,1)$, the \emph{$n$-th partial sum} of its 
Walsh--Fourier series is
\[
S_nf(x):=\sum_{k=0}^{n-1}\hat{f}(k)\,w_k(x)
=\int_0^1 f(t)\,D_n(x\dotplus t)\,dt,
\]
where $\hat{f}(k):=\int_0^1 f(t)\,w_k(t)\,dt$ is the $k$-th Walsh--Fourier coefficient, 
and $D_n(x):=\sum_{k=0}^{n-1}w_k(x)$ is the \emph{$n$-th Walsh--Dirichlet kernel}.
In particular, for any $n\in\mathbb{N}$, we have
\begin{equation}\label{eq:D2n}
D_{2^n}(x)=\begin{cases}
2^n&\text{if $x\in I_n$,}\\
0&\text{otherwise}.
\end{cases}
\end{equation}

For every $n\in\mathbb{P}$, define the \emph{$n$-th Walsh--Fej\'{e}r mean} of a 
function $f\in L^1[0,1)$ as
\[
\sigma_n f(x):=\frac{1}{n}\sum_{k=1}^{n} S_kf(x)=\int_0^1 f(t)\,K_n(x\dotplus t)\,dt,
\]
where the \emph{$n$-th Walsh--Fej\'{e}r kernel}---or simply 
\emph{Fej\'{e}r kernel}---is given by
\[
K_{n}(x):=\frac{1}{n}\sum_{k=1}^{n}D_{k}(x).
\]

The dyadic Fejér kernels $K_{2^n}$ exhibit much better properties than the general 
kernels~$K_n$. For any $n\in\mathbb{N}$ and $x\in[0,1)$, we have the exact 
formula~\cite[Lemma~3.2.1]{PTW22}
\begin{equation}\label{eq:K2n}
K_{2^n}(x)=
\begin{cases}
2^{n-1}+\frac12 & \textnormal{if }x\in I_n, \\
2^{i-1} & \textnormal{if } x\in I_n\dotplus 2^{-i-1} \textnormal{ for } 0\le i<n,\\
0 & \textnormal{otherwise},
\end{cases}
\end{equation}
which shows that $K_{2^n}\ge0$ everywhere on $[0,1)$. Comparing~\eqref{eq:D2n} to 
\eqref{eq:K2n} directly yields the estimate
\begin{equation}\label{ineq:DK}
D_{2^n}(x)\le 2\,K_{2^n}(x)
\end{equation}
for all $n\in\mathbb{N}$ and $x\in[0,1)$. Moreover, a classical identity
in~\cite[Chapter~1, Theorem~16(\emph{iii})]{SWS90} establishes the pointwise upper bound
\begin{equation}\label{est:K2n}
K_{2^{n}}(x)\le\sum_{i=0}^n 2^{i-n}D_{2^{n}}(x\dotplus2^{-i-1}).
\end{equation}

In general, for any $n\in\mathbb{P}$, the kernel $K_n$ 
satisfies~\cite[Corollary~3.2.7]{PTW22} the pointwise bound 
\begin{equation}\label{est:Kn}
|K_n(x)|\le \frac{C}{n}\sum_{i=0}^{|n|}2^i\,K_{2^i}(x),
\end{equation}
where $C>0$ is an absolute constant, and $|n|$ denotes the unique integer such that 
$2^{|n|}\le n<2^{|n|+1}$.

\section{Approximation in Weighted Lebesgue Spaces}

Hereafter, an integrable function $w$ on $[0,1)$ is called a \emph{weight} if
$0<w(x)<\infty$ almost everywhere. Recall that $\mathcal{D}$ denotes the collection of 
all dyadic intervals on $[0,1)$.
\begin{definition}\label{def:Muckenhoupt}
For $1<p<\infty$, a weight $w$ satisfies the $A_p(\mathcal{D})$ condition, 
written $w\in A_p(\mathcal{D})$, if
\[
[w]_{A_p}:=\sup_{I\in \mathcal{D}}
\bigg(\frac{1}{|I|}\int_{I}w(x)\,dx\bigg)
\bigg(\frac{1}{|I|}\int_{I}w(x)^{-\frac{1}{p-1}}\,dx\bigg)^{p-1}<\infty.
\]
For $p=1$, the weight $w$ satisfies the $A_1(\mathcal{D})$ condition, 
written $w\in A_1(\mathcal{D})$, if
\[
[w]_{A_1}:=\sup_{I\in\mathcal{D}}\,
\frac{\frac{1}{|I|}\int_I w(x)\,dx}{\essinf\limits_{x\in I}w(x)}
<\infty.
\]
\end{definition}
The $A_p(\mathcal{D})$ weights are the dyadic counterparts to the classical 
Muckenhoupt weights and share their foundational properties. For a measurable set 
$E\subseteq[0,1)$, denote 
\[
w(E):=\int_E w(x)\,dx.
\]
\begin{lemma}\label{le:A-inf}
The following properties hold:
\begin{enumerate}
\item[\textnormal{(i)}] A weight $w$ belongs to $\bigcup_{p\ge 1}A_p(\mathcal{D})$ if 
and only if it satisfies the $A_\infty(\mathcal{D})$ condition: there exist constants 
$C,\delta>0$ such that for any $I\in\mathcal{D}$ and measurable $E\subseteq I$, 
\[
\frac{w(E)}{w(I)}\le C\,\bigg(\frac{|E|}{|I|}\bigg)^{\delta}.
\]

\item[\textnormal{(ii)}] If $w\in A_1(\mathcal{D})$, there exists a constant $C>0$ 
such that for any $I\in\mathcal{D}$ and measurable $E\subseteq I$ with $|E|>0$,
\[
\frac{w(I)}{|I|}\le C\,\frac{w(E)}{|E|}.
\]
\end{enumerate}
\end{lemma}
\begin{proof}
The proof of~(i) is completely analogous to the classical proof by Coifman and 
Fefferman~\cite{CF74}, as detailed in~\cite[Chapter~IV, Lemma~2.5 and 
Theorem~2.9]{GCRF85}. In the context of Vilenkin systems, this property appeared 
in~\cite[p.~350]{JG75} and~\cite[Section~2]{Yo93}.

For~(ii), observe by definition of the $A_1(\mathcal{D})$ weights that
\[
\frac{w(I)}{|I|}\le C\,\essinf_{x\in I}w(x)
\le C\,\frac1{|E|}\int_E w(x)\,dx
=C\,\frac{w(E)}{|E|}
\]
which yields the desired inequality.
\end{proof}
Back in 1932, Paley~\cite{RP32} established that the Walsh system constitutes a basis in
$L^p[0,1)$ for $1<p<\infty$, relying on two fundamental bounds, now 
classically known as Paley's inequalities.
In~1974, Gundy and Wheeden~\cite[Theorem~2]{GW74} extended these inequalities to the 
weighted setting, laying the groundwork for the following result by 
Kazarian~\cite[Theorem~4, $(\alpha)\Leftrightarrow(\delta)$]{KK82}:
\begin{theorem}\label{th:Kazarian}
For a weight $w$ and $1<p<\infty$, the following are 
equivalent:
\begin{enumerate}
\item[\textnormal{(i)}] $w\in A_{p}(\mathcal{D})$.

\item[\textnormal{(ii)}] There exists a constant $C>0$, depending only on $p$ and 
$[w]_{A_{p}}$, such that for all $n\in\mathbb{N}$ and $f\in L^{p}_{w}[0,1)$,
\[
\int_0^1|S_{n}f(x)|^{p}\,w(x)\,dx
\le C\int_0^1 |f(x)|^p\,w(x)\,dx.
\]

\item[\textnormal{(iii)}] For all $f\in L^{p}_{w}[0,1)$,
\[
\lim_{n\to\infty}\int_0^1|S_nf(x)-f(x)|^p\,w(x)\,dx=0.
\]
\end{enumerate}	
\end{theorem}

Combining Theorem~\ref{th:Kazarian} with the pointwise inequality~\eqref{ineq:DK} 
between the Dirichlet and Fej\'er kernels yields the following equivalences for 
the Walsh--Fej\'er means:
\begin{theorem}\label{th:Fejer-Ap}
For a weight~$w$ and $1<p<\infty$, the following are equivalent:
\begin{enumerate}
\item[\textnormal{(i)}] $w\in A_{p}(\mathcal{D})$.

\item[\textnormal{(ii)}] There exists a constant $C>0$, depending only 
on~$p$ and $[w]_{A_p}$, such that for all $n\in\mathbb{P}$ and $f\in L^{p}_{w}[0,1)$,
\[
\int_0^1|\sigma_{n}f(x)|^p\,w(x)\,dx
\le C\int_0^1 |f(x)|^p\,w(x)\,dx.
\]

\item[\textnormal{(iii)}] For all $f\in L^p_w[0,1)$,
\[
\lim_{n\to\infty}\int_0^1|\sigma_n f(x)-f(x)|^p\,w(x)\,dx=0.
\]
\end{enumerate}
\end{theorem}
\begin{remark}
For general Vilenkin systems, results analogous to Theorem~\ref{th:Kazarian} were 
obtained by Gosselin~\cite{JG75} for the bounded case, and by 
Young~\cite[Theorem~1.1]{Yo93} for the unbounded case. Implication 
(i)$\,\Rightarrow\,$(ii) of Theorem~\ref{th:Fejer-Ap} can be found in 
\cite[Corollary~4.2]{FW25}; see also Weisz, Xie, and Yang~\cite{WXY26} for an extended 
study.
\end{remark}

The analogous result for the boundary case $p=1$ was previously unknown. 
As one of our main contributions, we prove the following theorem:
\begin{theorem}\label{th:Fejer-A1}
For a weight~$w$, the following are equivalent:
\begin{enumerate}
\item[\textnormal{(i)}] $w\in A_1(\mathcal{D})$.
	
\item[\textnormal{(ii)}] 
There exists a constant $C>0$, depending only on $[w]_{A_{1}}$, such that for all 
$n\in\mathbb{P}$ and $f\in L^1_w[0,1)$,
\[
\int_0^1|\sigma_{n}f(x)|\,w(x)\,dx
\le C\int_0^1 |f(x)|\,w(x)\,dx.
\]
	
\item[\textnormal{(iii)}] For all $f\in L^1_w [0,1)$,
\[
\lim_{n\to\infty}\int_0^1 |\sigma_nf(x)-f(x)|\,w(x)\,dx=0.
\]
\end{enumerate}
\end{theorem}
\begin{proof}
We first prove that (i)$\,\Rightarrow\,$(ii). Denote the dyadic intervals by
\[
I_{n,k}:=\bigg[\frac{k}{2^n},\frac{k+1}{2^n}\bigg)
\quad\textnormal{for } n,k\in\mathbb{N} 
\textnormal{ with }
0\le k<2^n.
\]
Notice that all points $x\in I_{n,k}$ share the same first $n$ dyadic coordinates
$x_0,x_1,\ldots,x_{n-1}$. Motivated by the Fej\'er kernel estimates~\eqref{est:K2n} 
and~\eqref{est:Kn}, we consider the operators
\[
U_n f(x)=\sum_{k=0}^{2^n-1}\chi_{I_{n,k}}(x)
\sum_{i=0}^{n} 2^{i-n}\frac{1}{|I^i_{n,k}|}
\int_{I^i_{n,k}}|f(t)|\,dt,
\]
where $I_{n,k}^i:=I_{n,k}\dotplus 2^{-i-1}$. 

The \emph{dyadic shift} $x\mapsto x\dotplus2^{-i-1}$ flips precisely the $i$-th 
coordinate of~$x$. For $0\le i<n$, this maps $I_{n,k}$ 
onto a distinct dyadic interval of the same length. For $i=n$, it maps $I_{n,k}$
onto itself, so that $I_{n,k}^n=I_{n,k}$. Because the shift is a bijection on $[0,1)$, 
it simply reorders the intervals. 
Therefore, for all $i\in\mathbb{N}$ satisfying $0\le i\le n$, we have
$|I^i_{n,k}|=|I_{n,k}|= 2^{-n}$ and 
\[
\bigcup_{k=0}^{2^n-1} I_{n,k}^i=[0,1).
\]

Let $n\in\mathbb{N}$, $f\in L^1_w[0,1)$, and $x\in[0,1)$. Applying 
estimate~\eqref{est:K2n} and the explicit formula~\eqref{eq:D2n} for the Dirichlet 
kernel, we have
\begin{align*}
\sigma_{2^n}|f(x)|&=\int_0^1|f(t)|\,K_{2^n}(x\dotplus t)\,dt
\\
&\le \sum_{i=0}^n 2^{i-n}
\int_0^1 |f(t)|\,D_{2^n}(x\dotplus t\dotplus 2^{-i-1})\,dt
\\
&=\sum_{i=0}^n 2^{i-n}
\int_0^1 |f(t\dotplus 2^{-i-1})|\,D_{2^n}(x\dotplus t)\,dt
\\
&=\sum_{k=0}^{2^n-1} \sum_{i=0}^n 2^{i-n}
\int_{I_{n,k}} |f(t\dotplus 2^{-i-1})|\,D_{2^n}(x\dotplus t)\,dt
\\
&=\sum_{k=0}^{2^n-1} \sum_{i=0}^n 2^{i-n} 
\int_{I_{n,k}} |f(t\dotplus 2^{-i-1})|\,2^n\,\chi_{I_{n,k}}(x)\,dt
\\
&=\sum_{k=0}^{2^n-1} \chi_{I_{n,k}}(x)
\sum_{i=0}^n 2^{i-n}\,2^n
\int_{I_{n,k}} |f(t\dotplus 2^{-i-1})|\,dt
\\
&=\sum_{k=0}^{2^n-1} \chi_{I_{n,k}}(x)
\sum_{i=0}^{n} 2^{i-n}\,2^n
\int_{I^i_{n,k}} |f(t)|\,dt
\\
&=\sum_{k=0}^{2^n-1} \chi_{I_{n,k}}(x)
\sum_{i=0}^n 2^{i-n}\,\frac{1}{|I^i_{n,k}|}
\int_{I^i_{n,k}} |f(t)|\,dt
=U_{n}f(x).
\end{align*}
Using~\eqref{est:Kn} to control the general Fej\'er means, for any $n\in\mathbb{P}$ we 
find
\[
|\sigma_n f(x)|\le 
\frac{C}{n}\sum_{i=0}^{|n|} 2^i\,\sigma_{2^i}|f(x)|
\le \frac{C}{n}\sum_{i=0}^{|n|}2^i\,U_{i}f(x),
\]
and consequently,
\begin{equation}\label{est:sigma-U}
\int_0^1 |\sigma_n f(x)|\,w(x)\,dx
\le \frac{C}{n}\sum_{i=0}^{|n|}2^i\int_0^1 U_i f(x)\,w(x)\,dx.
\end{equation}

It only remains to prove that for some constant $C>0$, independent of~$n\in\mathbb{N}$
and $f\in L^1_w[0,1)$, the following estimate holds:
\begin{equation}\label{est:U}
\int_0^1 U_n f(x)\,w(x)\,dx
\le C\int_0^1 |f(x)|\,w(x)\,dx.
\end{equation}
Once this is established, the desired inequality in~(ii) follows immediately 
from~\eqref{est:sigma-U}, since the geometric sum 
$\sum_{i=0}^{|n|}2^i=2^{|n|+1}-1$ is bounded by $2n$.

We now establish~\eqref{est:U}. Substituting the definition of~$U_nf$ and 
decomposing~$[0,1)$ into dyadic intervals, we note that for any fixed~$i$, the mapping 
$I_{n,k}\mapsto I_{n,k}^i$ simply reorders the intervals of length~$2^{-n}$. 
Re-indexing the sum over~$k$ thus seamlessly swaps the paired domains of integration:
\begin{align*}
\int_0^1 U_nf(x)\,w(x)\,dx
&=\sum_{k=0}^{2^n-1} \int_{I_{n,k}}w(x)\,dx\,
\sum_{i=0}^{n} 2^{i-n}
\frac{1}{|I_{n,k}^i|}\int_{I_{n,k}^i}|f(t)|\,dt
\\
&=\sum_{i=0}^n 2^{i-n} \sum_{k=0}^{2^n-1} 
\bigg(\frac{1}{|I_{n,k}^i|} \int_{I_{n,k}^i} |f(t)|\, dt\bigg)
\int_{I_{n,k}} w(x)\,dx
\\
&=\sum_{i=0}^{n} 2^{i-n} \sum_{k=0}^{2^n-1} 
\bigg(\frac{1}{|I_{n,k}|}\int_{I_{n,k}}|f(t)|\,dt\bigg)
\int_{I_{n,k}^i}w(x)\,dx
\\
&=
\sum_{k=0}^{2^n-1} 
\bigg(\frac{1}{|I_{n,k}|}\int_{I_{n,k}}|f(t)|\,dt\bigg)
\sum_{i=0}^n 2^{i-n}\,w(I_{n,k}^i).
\end{align*}

Since the dyadic shift by~$2^{-i-1}$ alters only the $i$-th coordinate,
the intervals $I_{n,k}$ and $I_{n,k}^i$ share the first $i$ coordinates. Consequently,
for any~$k$, there exists a unique dyadic interval $I_{i,l}$ of length~$2^{-i}$ that 
contains both $I_{n,k}$ and $I_{n,k}^i$. Then Lemma~\ref{le:A-inf}, applied to 
inclusions $I_{n,k}^i\subseteq I_{i,l}$ and $I_{n,k}\subseteq I_{i,l}$, yields
\begin{align*}
w(I_{n,k}^i) 
&\le C\,\bigg(\frac{|I_{n,k}^i|}{|I_{i,l}|}\bigg)^\delta\,w(I_{i,l})
\le C\,2^{(i-n)\delta}\,\frac{|I_{i,l}|}{|I_{n,k}|}\,w(I_{n,k})
\\
&\le C\,2^{(i-n)(\delta-1)}\,w(I_{n,k}).
\end{align*}
For any $k\in\mathbb{N}$ satisfying $0\le k<2^n$, this results in
\begin{align}\label{eq:sum}
\sum_{i=0}^n 2^{i-n}\,w(I_{n,k}^i)
&\le C\sum_{i=0}^n 2^{i-n}\, 2^{(i-n)(\delta-1)}\,w(I_{n,k})
\nonumber
\\
&\le C\,\sum_{j=0}^\infty 2^{-j\delta}\,w(I_{n,k})
\le C\,w(I_{n,k}).
\end{align}
Combining~\eqref{eq:sum} with the $A_1(\mathcal{D})$ condition for the weight~$w$,
we finally obtain
\begin{align*}
\int_0^1 U_n f(x)\,w(x)\,dx
&\le \sum_{k=0}^{2^n-1}
\frac{\frac{1}{|I_{n,k}|} \int_{I_{n,k}}|f(t)|\,w(t)\,dt}
{\essinf\limits_{x\in I_{n,k}}w(x)}\,
\sum_{i=0}^n 2^{i-n}\,w(I_{n,k}^i)
\\
&\le C\,\sum_{k=0}^{2^n-1}
\int_{I_{n,k}}|f(t)|\,w(t)\,dt
\cdot 
\frac{\frac1{|I_{n,k}|}\,w(I_{n,k})}{\essinf\limits_{x\in I_{n,k}}w(x)}
\\
&\le C\,\sum_{k=0}^{2^n-1}
\int_{I_{n,k}}|f(t)|\,w(t)\,dt
=C\int_0^1 |f(x)|\,w(x)\,dx,
\end{align*}
thus establishing~\eqref{est:U}. The proof of (i)$\,\Rightarrow\,$(ii) is complete.

Let us prove the reverse implication (ii)$\,\Rightarrow\,$(i). 
In view of inequality~\eqref{ineq:DK}, hypothesis~(ii) implies that 
for all $n\in\mathbb{N}$ and all nonnegative functions $f\in L^1_w[0,1)$, 
there holds
\begin{equation}\label{ineq:S2n}
\int_0^1 S_{2^n}f(x)\,w(x)\,dx 
\le C \int_0^1 f(x)\,w(x)\,dx.
\end{equation}
Recall that the Dirichlet kernel $D_{2^n}(x\dotplus t)$ equals $2^n$ when $x$ and $t$ 
belong to the same dyadic interval $I_{n,k}$, and $0$ otherwise. Thus, the dyadic 
partial sum is simply a step function:
\begin{equation}\label{eq:S2n}
S_{2^n}f(x)=\int_0^1 f(t)\,D_{2^n}(x\dotplus t)\,dt
=\sum_{k=0}^{2^n-1} \chi_{I_{n,k}}(x) 
\frac{1}{|I_{n,k}|} \int_{I_{n,k}} f(t)\,dt.
\end{equation}
Substituting this explicit formula into the left-hand side of~\eqref{ineq:S2n} yields
\begin{align*}
\int_0^1 S_{2^n}f(x)\,w(x)\,dx&=
\sum_{k=0}^{2^n-1}
\bigg(\frac{1}{|I_{n,k}|} \int_{I_{n,k}} f(t)\,dt\bigg)
\int_{I_{n,k}} w(x)\,dx
\\
&=\sum_{k=0}^{2^n-1}
\int_{I_{n,k}}f(t)\,\frac{w(I_{n,k})}{|I_{n,k}|}\,dt
\\
&=\int_0^1 f(t)\,
\bigg(\sum_{k=0}^{2^n-1}\chi_{I_{n,k}}(t)\,\frac{w(I_{n,k})}{|I_{n,k}|}\bigg)\,dt.
\end{align*}
In order for inequality~\eqref{ineq:S2n} to hold for all $f\ge0$, the multiplier in the 
parentheses must be pointwise bounded by $C\,w(t)$ for almost every $t\in[0,1)$.
Restricting $t$ to an arbitrary dyadic interval $I_{n,k}$ gives
\[
\frac{w(I_{n,k})}{|I_{n,k}|}\le C\,w(t)
\quad\textnormal{for almost every }t\in I_{n,k}. 
\]
Taking the essential infimum over $t\in I_{n,k}$ on the right-hand side 
establishes that $w\in A_1(\mathcal{D})$.

It remains to show the equivalence (ii)$\,\Leftrightarrow\,$(iii). Simple functions 
are dense in $L^1_w[0,1)$~\cite[Chapter~4, Lemma~3.4]{BS88}. 
Since $w \in L^1[0,1)$, absolute continuity of the weighted measure ensures that 
any simple function can be approximated in $L^1_w$-norm by a function that is constant 
on dyadic intervals of length $2^{-n}$ for some sufficiently large $n$. 
Consequently, the dyadic step functions
\[
g(x)=\sum_{k=0}^{2^n-1}c_k\,\chi_{I_{n,k}}(x)=
\sum_{j=0}^{2^n-1} a_j\,w_j(x), 
\]
which are precisely the Walsh polynomials of degree less than $2^n$, form a 
dense subspace of~$L^1_w[0,1)$.

Assume (ii), meaning there is a constant $C>0$ such that 
$\|\sigma_nf\|_{L^1_w}\le C\,\|f\|_{L^1_w}$ for all $n\in\mathbb{P}$. 
Fix $f\in L^1_w[0,1)$ and $\varepsilon>0$. Choose a Walsh polynomial~$g$ such that 
$\|f-g\|_{L^1_w}<\varepsilon$. Since $S_ng=g$ for all sufficiently large $n$, 
the regularity of the Fej\'er means ensures that $\sigma_n g\to g$ in~$L^1_w$.
Thus, $\|\sigma_ng-g\|_{L^1_w}<\varepsilon$ for large $n$.
The convergence~(iii) then follows directly from the inequality
\[
\|\sigma_nf-f\|_{L^1_w}\le 
\|\sigma_n(f-g)\|_{L^1_w}
+\|\sigma_ng-g\|_{L^1_w}
+\|f-g\|_{L^1_w}
<(C+2)\,\varepsilon.
\]

Conversely, assume~(iii). Since every convergent sequence in a normed space is bounded,
for each fixed $f\in L^1_w[0,1)$ we have 
\[
\sup_{n\in\mathbb{P}}\|\sigma_nf\|_{L^1_w}<\infty. 
\]
Since $L^1_w[0,1)$ is a Banach space, the uniform boundedness principle guarantees that 
the operators $\sigma_n$ are uniformly bounded on $L^1_w[0,1)$, which is exactly~(ii).
\end{proof}

\section{Extrapolation to Variable Lebesgue Spaces}
\label{sec:var-Lebesgue}

We recall the definition of variable Lebesgue spaces. 
A \emph{variable exponent} $p(\cdot):[0,1)\to[1,\infty]$ is a measurable function, 
whose range is described by the essential bounds
\[
p_-:=\essinf_{x\in[0,1)}p(x)
\quad\text{and}\quad
p_+:=\esssup_{x\in[0,1)}p(x).
\]
The conjugate exponent $p'(\cdot)$ is defined by the pointwise identity
\[
\frac1{p(x)}+\frac1{p'(x)}=1,
\quad x\in[0,1),
\]
with the standard convention that $1/\infty=0$.

Let $t^\infty:=0$ for $t\in[0,1]$, and $t^\infty:=\infty$ for $t>1$. Given an
exponent~$p(\cdot)$, the \emph{variable Lebesgue space} $L^{p(\cdot)}[0,1)$ consists of
all measurable functions~$f$ on $[0,1)$ such that for some $\lambda>0$,
\[
\varrho_{p(\cdot)}(f/\lambda):=
\int_0^1 \bigg|\frac{f(x)}{\lambda}\bigg|^{p(x)} dx<\infty.
\]
This forms a Banach space when equipped with the Luxemburg norm
\[
\|f\|_{p(\cdot)}:=\inf
\big\{\lambda>0\ :\ \varrho_{p(\cdot)}(f/\lambda)\le1 \big\}.
\]

The following condition on the exponent~$p(\cdot)$ generalizes
the classical Muckenhoupt $A_p(\mathcal{D})$ condition from 
Definition~\ref{def:Muckenhoupt}:
\begin{definition}
A variable exponent~$p(\cdot)$ satisfies the $\mathcal{A}(\mathcal{D})$ condition, 
written $p(\cdot)\in\mathcal{A}(\mathcal{D})$, if 
\[
[p]_{\mathcal{A}}:=\sup_{I\in\mathcal{D}}
\frac{1}{|I|}\,\|\chi_I\|_{p(\cdot)} \|\chi_I\|_{p'(\cdot)}<\infty.
\]
\end{definition}
This condition completely characterizes the boundedness of the 
\emph{dyadic maximal operator}~$M^\mathcal{D}$ on the spaces 
$L^{p(\cdot)}[0,1)$. For $f\in L^1[0,1)$, the operator is defined by 
\[
M^\mathcal{D}f(x):=\sup_{\substack{I\ni x,\\ I\in\mathcal{D}}}
\frac{1}{|I|} \int_I |f(t)|\,dt,
\quad x\in[0,1),
\]
where the supremum is taken over all dyadic intervals~$I$ on~$[0,1)$ containing~$x$.
The fundamental characterization is given by the following theorem, which extends
the classical result by Kopaliani~\cite[Theorem~1.2]{TK07}:
\begin{theorem}\label{th:Kopaliani}
Let an exponent $p(\cdot)$ satisfy $1<p_-\le p_+\le\infty$. The dyadic 
maximal operator~$M^\mathcal{D}$ is bounded in $L^{p(\cdot)}[0,1)$ 
\emph{if and only if} $p(\cdot)\in\mathcal{A}(\mathcal{D})$.
\end{theorem}
\begin{proof}
The necessity of the $\mathcal{A}(\mathcal{D})$ condition follows exactly as in the 
proof of~\cite[Theorem~4.5.7, (a)$\,\Rightarrow\,$(c)]{DHHR11} using the norm 
conjugate formula. 

For the sufficiency, we adapt the approach 
from~\cite[Theorem~1.1 and Proposition~5.4]{ADK26}, which accommodates the case
$p_+\le\infty$. The proof simplifies in our setting: the Nekvinda decay 
condition is trivially satisfied on the bounded domain~$[0,1)$, and analyzing 
$M^\mathcal{D}$ directly eliminates the need for shifted grids, requiring 
only~$\mathcal{A}(\mathcal{D})$. The boundedness of the maximal operator over the 
single grid $\mathcal{D}$ is first established for $p_+<\infty$, and then extended to 
$p_+\le\infty$ via pointwise exponent approximation and Fatou's lemma.
\end{proof}
\begin{remark}
The behavior of maximal operators with respect to general dyadic systems in 
variable Lebesgue and Musielak--Orlicz spaces has been the subject of active recent 
research; see, e.g., Shalukhina~\cite{AS26}.
\end{remark}

A powerful tool for deducing norm estimates on variable Lebesgue spaces from 
classical weighted inequalities is the Rubio de Francia extrapolation theorem:
\begin{theorem}\label{th:RdF}
Let~$\mathcal{F}$ be a family of pairs of nonnegative, measurable functions. Suppose 
that for some $p_0\ge1$ and all weights $w\in A_1(\mathcal{D})$, 
\begin{equation}\label{eq:RdF-hyp}
\int_0^1 f(x)^{p_0}\,w(x)\,dx
\le C\int_0^1 g(x)^{p_0}\,w(x)\,dx
\end{equation}
holds for all $(f,g)\in\mathcal{F}$ whenever the left-hand side is finite, where $C$ 
depends only on~$p_0$ and $[w]_{A_1}$. If an exponent $p(\cdot)$ satisfies 
$p_0\le p_-\le p_+<\infty$ and $M^\mathcal{D}$ is bounded 
on~$L^{(p(\cdot)/p_0)'}[0,1)$, then 
\[
\|f\|_{p(\cdot)}\le C_{p(\cdot)}\,\|g\|_{p(\cdot)}
\]
for all $(f,g)\in\mathcal{F}$ with $f\in L^{p(\cdot)}[0,1)$, where 
$C_{p(\cdot)}$ depends only on~$p(\cdot)$.
\end{theorem}
\begin{remark}
The extrapolation theorem for variable Lebesgue spaces was originally 
established in~\cite[Theorem~1.3]{CUFMP06} for the standard maximal 
operator over Euclidean cubes. It was later extended in~\cite[Theorem~4.6]{CUMP11} 
to general Banach function spaces and maximal operators defined 
via arbitrary \emph{Muckenhoupt bases}~\cite[Definition~3.1]{CUMP11}.
Since the collection of dyadic intervals~$\mathcal{D}$ forms such a 
basis, Theorem~\ref{th:RdF} follows directly from this generalized result
setting $\mathbb{X}:=L^{p(\cdot)}[0,1)$.
\end{remark}

The Rubio de Francia extrapolation allows us to extend 
Theorems~\ref{th:Fejer-Ap} and~\ref{th:Fejer-A1} to variable Lebesgue spaces as follows.
\begin{theorem}\label{th:main-variable}
For an exponent $p(\cdot)$ with $1\le p_-\le p_+<\infty$, 
the following are equivalent:
\begin{enumerate}
\item[\textnormal{(i)}] $p(\cdot)\in\mathcal{A}(\mathcal{D})$.

\item[\textnormal{(ii)}] There exists a constant $C>0$, depending only on $p(\cdot)$
and $[p]_\mathcal{A}$, such that for all $n\in\mathbb{P}$ and $f\in L^{p(\cdot)}[0,1)$,
\[
\|\sigma_{n}f\|_{p(\cdot)}\le C\,\|f\|_{p(\cdot)}.
\]
	
\item[\textnormal{(iii)}] For all $f\in L^{p(\cdot)}[0,1)$,
\[
\lim_{n\to\infty}\|\sigma_{n}f-f\|_{p(\cdot)}=0.
\]
\end{enumerate}
\end{theorem}
\begin{proof}
We first establish (i)$\,\Rightarrow\,$(ii). Because $[0,1)$ has 
finite measure and $p_-\ge1$, the embedding 
$L^{p(\cdot)}[0,1)\subseteq L^1[0,1)$~\cite[Corollary~3.3.4]{DHHR11} 
ensures that the Fej\'er means $\sigma_n f$ are well-defined.

Let $\mathcal{F}$ be the family of pairs
$(|\sigma_nf|,|f|)$ for all $n\in\mathbb{P}$ and $f\in L^{p(\cdot)}[0,1)$.
By Theorem~\ref{th:Fejer-A1}, this family satisfies the $A_1(\mathcal{D})$ weighted 
inequality required by~\eqref{eq:RdF-hyp} with $p_0=1$. Furthermore, 
assumption~(i) implies that $p'(\cdot)\in\mathcal{A}(\mathcal{D})$. Since  
$1<(p')_-\le (p')_+\le\infty$, Theorem~\ref{th:Kopaliani} ensures that 
$M^\mathcal{D}$ is bounded on $L^{p'(\cdot)}[0,1)$. Applying Theorem~\ref{th:RdF}
then immediately yields~(ii).

To prove (ii)$\,\Rightarrow\,$(i), we argue as in Theorem~\ref{th:Fejer-A1}.  
Pointwise domination of the dyadic partial sums by the Fejér means applies to 
non-negative functions. Therefore, hypothesis~(ii) implies that there exists
$C>0$ such that 
\[
\|S_{2^n}f\|_{p(\cdot)}\le C\,\|f\|_{p(\cdot)}
\]
for all $f\ge 0$ in $L^{p(\cdot)}[0,1)$. Using the norm conjugate 
formula~\cite[Corollary~3.2.14]{DHHR11} and the explicit step-function 
expression~\eqref{eq:S2n} for $S_{2^n}f$, for any dyadic interval $I_{n,k}$ we obtain:
\begin{align*}
\frac1{|I_{n,k}|}\,\|\chi_{I_{n,k}}\|_{p(\cdot)} \|\chi_{I_{n,k}}\|_{p'(\cdot)}
&\le 2\sup_{\|f\|_{p(\cdot)}\le1}
\bigg\| 
\bigg(\frac{1}{|I_{n,k}|}\int_{I_{n,k}}|f(t)|\,dt\bigg) \chi_{I_{n,k}}
\bigg\|_{p(\cdot)}
\\
&\le 2\sup_{\|f\|_{p(\cdot)}\le1}\|S_{2^n}|f|\,\|_{p(\cdot)}\le 2C.
\end{align*}
Taking the supremum over all intervals $I_{n,k}\in\mathcal{D}$ yields 
$p(\cdot)\in\mathcal{A}(\mathcal{D})$, which is exactly~(i).

Finally, the equivalence (ii)$\,\Leftrightarrow\,$(iii) follows
as in Theorem~\ref{th:Fejer-A1}. The required density of Walsh polynomials 
in~$L^{p(\cdot)}[0,1)$ holds because $p_+<\infty$ 
ensures the density of simple functions~\cite[Corollary~3.4.10]{DHHR11}.
\end{proof}
\begin{remark}
Adamadze and Kopaliani~\cite[Theorem~1.10]{AK23} established an analogue of 
Theorem~\ref{th:main-variable} for the partial sums~$S_n$ of general Vilenkin--Fourier 
series. Their characterization is restricted to the range $1<p_-\le p_+<\infty$.
\end{remark}

\section{Walsh--Fej\'er Means in Variable Lebesgue Spaces}
\label{sec:Fejer}

A fundamental property of the Walsh--Fej\'er means~$\sigma_n$ is their uniform 
boundedness on $L^p[0,1)$ for all $1 \le p \le \infty$. By 
Young's convolution inequality, we have $\|\sigma_nf\|_p\le \|f\|_p\,\|K_n\|_1$. 
This and the uniform $L^1$-boundedness of the 
kernels~$K_n$~\cite[Chapter~1, Theorem~16(\emph{v})]{SWS90} guarantee the 
existence of an absolute constant $C>0$ such that
\[
\|\sigma_n f\|_p \le C\,\|f\|_p.
\]

Recall that the dyadic \emph{modulus of continuity} of a function $f\in L^p[0,1)$ is 
defined by
\[
\omega_p(f,2^{-k}):= \sup_{0\le h<2^{-k}}
\| f(\cdot \dotplus h) - f(\cdot) \|_p,
\quad k\in\mathbb{N}.
\]
For $1\le p<\infty$, it was established in~\cite[Theorem~3.3.1]{PTW22} that
for all $n\in\mathbb{P}$ and $f\in L^p[0,1)$,
\begin{equation}\label{ineq:Fej-Lp}
\|\sigma_n f-f\|_{p}
\le C\sum_{k=0}^{|n|} 2^{k-|n|}\,\omega_p (f,2^{-k}),
\end{equation}
where $C>0$ depends only on 
$p$ and, as before, $|n|$ denotes the unique integer satisfying 
$2^{|n|}\le n<2^{|n|+1}$. As a consequence of~\eqref{ineq:Fej-Lp}, 
if a function $f$ belongs to the Lipschitz class
\[
\textnormal{Lip}(\alpha,p):=
\big\{f\in L^p[0,1)\,:\, 
\omega_p (f,2^{-k})=O(2^{-\alpha k})
\textnormal{ as }k\to\infty \big\},
\]
for some $\alpha>0$, its approximation error 
satisfies~\cite[Chapter~5, Theorem~4]{SWS90}
\[
\|\sigma_nf-f\|_p=
\begin{cases}
O(n^{-1}) & \textnormal{if }\alpha>1, \\
O(n^{-1}\log n) & \textnormal{if }\alpha=1, \\
O(n^{-\alpha}) & \textnormal{if }0<\alpha<1
\end{cases}
\quad\textnormal{as }n\to\infty.
\]

In this section, we extend inequality~\eqref{ineq:Fej-Lp} and its asymptotic corollary 
to variable Lebesgue spaces. In classical Lebesgue spaces, the dyadically 
shifted function $f(\cdot\dotplus h)$ always belongs to the same space as $f$. In the 
variable setting, however, this translation invariance fails: a shift 
$f(\cdot\dotplus h)$ might not belong to $L^{p(\cdot)}[0,1)$ even if $f$ does.
Consequently, the classical modulus of continuity is generally undefined.

To resolve this, we rely on the two-sided 
inequality from~\cite[Chapter~5, Theorem~2]{SWS90}, which guarantees that 
for $1\le p<\infty$,
\[
\frac12\,\omega_p(f,2^{-k})\le \| f-S_{2^k}f \|_p \le \omega_p (f,2^{-k}), 
\]
and define the \emph{modulus of smoothness} for $f\in L^{p(\cdot)}[0,1)$ directly 
via the approximation error of its dyadic partial sums:
\[
\omega_{p(\cdot)}(f,2^{-k}):=\|f-S_{2^k}f\|_{p(\cdot)},
\quad k\in\mathbb{N}.
\]

With this definition in place, we prove the following theorem.
\begin{theorem}\label{th:Fej-var-Lp}
Let an exponent $p(\cdot)\in\mathcal{A}(\mathcal{D})$ satisfy $1\le p_-\le p_+<\infty$.
Then there exists a constant $C>0$, depending only on $p(\cdot)$ and $[p]_\mathcal{A}$, 
such that for all $n\in\mathbb{P}$ and $f\in L^{p(\cdot)}[0,1)$,
\[
\|\sigma_n f-f\|_{p(\cdot)}
\le C\sum_{k=0}^{|n|} 2^{k-|n|}\,\omega_{p(\cdot)} (f,2^{-k}),
\]
where $|n|\in\mathbb{N}$ is given by $2^{|n|}\le n<2^{|n|+1}$.
\end{theorem}
\begin{proof}
For an arbitrary function $f\in L^{p(\cdot)}[0,1)$, decompose the error
of its approximation by the Fej\'er means $\sigma_n f$ into three parts:
\begin{align*}	
f-\sigma_nf &= \frac1n \sum_{i=1}^n (f-S_i f)
\\
&=\frac1n \bigg(\sum_{k=0}^{|n|-1} \sum_{i=2^k}^{2^{k+1}-1}
(f-S_i f)
+\sum_{i=2^{|n|}}^{n} (f-S_i f)\bigg)
\\
&=\frac1n \bigg(\sum_{k=0}^{|n|-1}\sum_{i=2^k}^{2^{k+1}-1}
(f-S_{2^k}f)
+\sum_{i=2^{|n|}}^{n} (f-S_{2^{|n|}}f) \bigg)
\\
&\quad -\frac1n \sum_{k=0}^{|n|-1} \sum_{i=2^k}^{2^{k+1}-1}
(S_if-S_{2^k}f)
-\frac1n \sum_{i=2^{|n|}}^{n} (S_if-S_{2^{|n|}}f)
\\	
&=:T_1-T_2-T_3,
\end{align*}
where the terms $T_1$, $T_2$, and $T_3$ correspond to the three resulting sums.

Evaluating the inner sums in $T_1$ yields
\[
T_1=
\frac1n \sum_{k=0}^{|n|-1} 2^k\,(f-S_{2^k}f) 
+\frac{n-2^{|n|}+1}{n}\,(f-S_{2^{|n|}}f). 
\]
Applying the triangle inequality, 
the bounds $2^{|n|}\le n<2^{|n|+1}$, and the definition of the modulus of 
smoothness, we obtain the first estimate
\begin{align}\label{est:T1}
\|T_1\|_{p(\cdot)}&\le
\sum_{k=0}^{|n|-1} \frac{2^k}{2^{|n|}}\,\|f-S_{2^k}f\|_{p(\cdot)} 
+\|f-S_{2^{|n|}}f\|_{p(\cdot)}
\nonumber
\\
&\le \sum_{k=0}^{|n|} 2^{k-|n|}\,\omega_{p(\cdot)}(f,2^{-k}).
\end{align}
	
To bound $T_2$, we observe that for $2^k\le i< 2^{k+1}$, the nesting property 
$S_i(S_{2^k}f) = S_{2^k}f$ of the partial sums gives the identity
\[
S_if-S_{2^k}f=S_if-S_i(S_{2^k}f)=S_i(f-S_{2^k}f).
\]
Rearranging the inner sum in terms of the Fej\'er means, we rewrite $T_2$ as
\begin{align*}
T_2 &= \frac1n \sum_{k=0}^{|n|-1} \sum_{i=2^k}^{2^{k+1}-1}
S_i (f-S_{2^k}f)
\\
&= \frac1n \sum_{k=0}^{|n|-1}
\big((2^{k+1}-1)\,\sigma_{2^{k+1}-1}(f-S_{2^k}f)
-(2^k-1)\,\sigma_{2^k-1}(f-S_{2^k}f)\big).
\end{align*}
Since $p(\cdot)\in\mathcal{A}(\mathcal{D})$ satisfies $1\le p_-\le p_+<\infty$, 
it follows from Theorem~\ref{th:main-variable} that the operators $\sigma_n$ are 
uniformly bounded on $L^{p(\cdot)}[0,1)$. Let $C_1>0$ be the uniform operator bound, 
depending only on $p(\cdot)$ and $[p]_\mathcal{A}$. Taking the norm of $T_2$ yields
\begin{align}\label{est:T2}	
\|T_2\|_{p(\cdot)}
&\le \frac{C_1}{n} \sum_{k=0}^{|n|-1}
\big(2^{k+1}\,\|f-S_{2^k}f\|_{p(\cdot)}
+2^k\,\|f-S_{2^k}f\|_{p(\cdot)}\big)
\nonumber
\\
&\le \frac{C_1}{2^{|n|}} \sum_{k=0}^{|n|-1}
(2^{k+1}+2^k)\,\omega_{p(\cdot)}(f,2^{-k})
\nonumber
\\
&= 3C_1 \sum_{k=0}^{|n|-1}2^{k-|n|}\,\omega_{p(\cdot)}(f,2^{-k}).
\end{align}

Similarly, we express $T_3$ via the Fej\'er means:
\begin{align*}
T_3&=\frac1n \sum_{i=2^{|n|}}^{n} S_i(f-S_{2^{|n|}}f)
\\
&=\frac1n \big(n\,\sigma_n(f-S_{2^{|n|}}f)
-(2^{|n|}-1)\,\sigma_{2^{|n|}-1}(f-S_{2^{|n|}}f)\big).
\end{align*}
Applying the uniform operator bound $C_1$ again gives
\begin{align}\label{est:T3}
\|T_3\|_{p(\cdot)}&\le 
C_1\,\|f-S_{2^{|n|}}f\|_{p(\cdot)}+
\frac{2^{|n|}-1}{n}\,C_1\,\|f-S_{2^{|n|}}f\|_{p(\cdot)}
\nonumber
\\
&\le 2C_1\,\omega_{p(\cdot)}(f,2^{-|n|}).
\end{align}

By the triangle inequality, $\|f-\sigma_n f\|_{p(\cdot)} \le 
\|T_1\|_{p(\cdot)}+\|T_2\|_{p(\cdot)}+\|T_3\|_{p(\cdot)}$.
Summing the three norm estimates~\eqref{est:T1}, \eqref{est:T2}, and
\eqref{est:T3} yields the desired inequality and completes the proof.
\end{proof}

\begin{corollary}\label{cor:Lip-Fejer}
Let an exponent $p(\cdot)\in\mathcal{A}(\mathcal{D})$ satisfy $1\le p_-\le p_+<\infty$.
If $f\in\textnormal{Lip}(\alpha,p(\cdot))$ for some $\alpha>0$, that is,
\[
\omega_{p(\cdot)}(f,2^{-k})=O(2^{-\alpha k})
\textnormal{ as } k\to\infty,
\]
then
\[
\|\sigma_nf-f\|_{p(\cdot)}=
\begin{cases}
O(n^{-1}) & \textnormal{if }\alpha>1, \\
O(n^{-1}\log n) & \textnormal{if }\alpha=1, \\
O(n^{-\alpha}) & \textnormal{if }0<\alpha<1
\end{cases}
\quad\textnormal{as }n\to\infty.
\]
\end{corollary}

\section{Application to Walsh--N\"orlund Means}

N\"orlund means generalize the Fej\'er means 
discussed in Section~\ref{sec:Fejer}. Let $\{q_k\}_{k\in\mathbb{N}}$ be a sequence of 
nonnegative numbers with $q_0>0$ whose partial sums 
\[
Q_n:=\sum_{k=0}^{n-1}q_k,
\quad n\in\mathbb{P},
\]
satisfy $\lim_{n\to\infty} Q_n=\infty$. For a function $f\in L^1[0,1)$, 
the \emph{$n$-th N\"orlund mean} of its Walsh--Fourier series is defined as
\begin{equation}\label{eq:tn}
t_nf(x):=\frac1{Q_n}\sum_{k=1}^n q_{n-k}\,S_kf(x)=\int_0^1 f(t)\,F_n(x\dotplus t)\,dt,
\end{equation}
where the \emph{n-th N\"orlund kernel} is given by
\[
F_n(x):=\frac1{Q_n} \sum_{k=1}^n q_{n-k}\,D_k(x).
\]
Setting $q_k = 1$ for all $k$ recovers the Fej\'er means $\sigma_n f$. 
The convergence of the N\"orlund means has been extensively studied; we refer the 
reader to the works of Fridli, Manchanda, and Siddiqi~\cite{FMS08}, Blahota and 
Nagy~\cite{BN18,BNT19}, Baramidze et al.~\cite{BNPT23, BPTT23}, and 
the monograph by Persson, Tephnadze, and Weisz~\cite{PTW22}.

For $t_n f$ to successfully approximate $f$, the N\"orlund summation method must 
be \emph{regular}, meaning that $t_n f$ converges to $f$ whenever the partial sums 
$S_k f$ do.
The method~\eqref{eq:tn}
is regular if and only if~\cite[Theorem~4.3.2]{PTW22} 
\[
\lim_{n\to\infty}\frac{q_{n-1}}{Q_n}=0.
\]
This condition is fundamentally tied to the monotonicity of~$\{q_k\}_{k\in\mathbb{N}}$. 
While every non-increasing sequence yields a regular 
method, non-decreasing sequences may not~\cite[Corollary~4.3.3]{PTW22}. 
We restrict our attention exclusively to regular N\"orlund means.

By applying the Abel transformation, the N\"orlund means can be expressed 
as a weighted linear combination of the Fej\'er means~\cite[Identity~(4.26)]{PTW22}. 
Taking the norm of this decomposition allows us to deduce the uniform boundedness 
of~$t_n$ directly from that of~$\sigma_n$, once $\{q_k\}_{k\in\mathbb{N}}$ satisfies 
the requisite monotonicity and regularity conditions. 
Combining Theorems~\ref{th:Fejer-Ap} 
and~\ref{th:Fejer-A1} thus gives the following result:
\begin{theorem}\label{th:Norlund-Ap}
Let~$w$ be a weight on~$[0,1)$ and $1\le p<\infty$. Suppose that $t_n$ is a regular 
N\"orlund mean generated by a sequence $\{q_k\}_{k\in\mathbb{N}}$ that is either 
non-decreasing and satisfies 
\begin{equation}\label{cond:non-decr}
\frac{q_{n-1}}{Q_n}=O\bigg(\frac{1}{n}\bigg)
\quad\textnormal{as }n\to\infty,
\end{equation} 
or non-increasing and satisfies
\begin{equation}\label{cond:non-incr}
\frac1{Q_n}=O\bigg(\frac1n\bigg)
\quad\textnormal{as }n\to\infty.
\end{equation} 
Then the following are equivalent:
\begin{enumerate}
\item[\textnormal{(i)}] $w\in A_{p}(\mathcal{D})$.

\item[\textnormal{(ii)}] There exists a constant $C>0$, depending only 
on~$p$ and $[w]_{A_p}$, such that 
for all $n\in\mathbb{P}$ and $f\in L^{p}_{w}[0,1)$,
\[
\int_0^1|t_n f(x)|^p\,w(x)\,dx
\le C\int_0^1 |f(x)|^p\,w(x)\,dx.
\]

\item[\textnormal{(iii)}] For all $f\in L^p_w[0,1)$,
\[
\lim_{n\to\infty}\int_0^1|t_n f(x)-f(x)|^p\,w(x)\,dx=0.
\]
\end{enumerate}
\end{theorem}
\begin{remark}
For non-decreasing sequences, condition~\eqref{cond:non-decr} is not required for the 
implication (i)$\,\Rightarrow\,$(ii). It is, however, crucial for the reverse 
implication (ii)$\,\Rightarrow\,$(i). This direction relies on the lower bound for the
N\"orlund kernel
\[
F_{2^n}(x\dotplus t)\ge c\,2^n
\quad\textnormal{for all }x,t\in I_{n,k},
\]
where $c>0$ depends only on the constant in~\eqref{cond:non-decr}, and $I_{n,k}$ is a 
dyadic interval of length~$2^{-n}$. 

Let us check the bound. For $x,t\in I_{n,k}$, their dyadic sum $y=x\dotplus t$ falls 
into~$I_n$. For $1\le k\le 2^n$, we have $D_k(y)=k$. In the non-decreasing case, 
assumption~\eqref{cond:non-decr} implies the existence of $l\in\mathbb{P}$ such that
for all sufficiently large~$n$, we have $2^{n-l}\,q_{2^n-1}\le Q_{2^n}/2$. 
Since $\{q_k\}_{k\in\mathbb{N}}$ is non-decreasing, the sum of any $r$ consecutive 
terms ending at $q_{2^n-1}$ is at most $r\,q_{2^n-1}$. Thus, for each $r\in\mathbb{N}$ 
satisfying $0\le r< 2^{n-l}$,
\[
Q_{2^n-r}\ge Q_{2^n}-r\,q_{2^n-1}\ge\frac{Q_{2^n}}{2}.
\]
Substituting $D_k(y)=k$ and using the above estimate, we obtain
\begin{align*}
F_{2^n}(y)&=\frac1{Q_{2^n}}\sum_{k=1}^{2^n} q_{2^n-k}\,k
=\frac1{Q_{2^n}}\sum_{m=1}^{2^n}Q_m
\\
&\ge\frac{1}{Q_{2^n}}\sum_{r=0}^{2^{n-l}-1} Q_{2^n-r}\ge\frac{2^n}{2^{l+1}}.
\end{align*}
Setting $c = 2^{-(l+1)}$ establishes the desired estimate $F_{2^n}(y)\ge c\,2^n$.

In the case of a non-increasing sequence $\{q_k\}_{k\in\mathbb{N}}$, 
the sum of the second half of the terms in $Q_{2^n}$ is bounded by the first half, 
yielding $Q_{2^n}\le 2\,Q_{2^{n-1}}$. We can therefore simply estimate
\[
F_{2^n}(y)=\frac{1}{Q_{2^n}}\sum_{m=1}^{2^n} Q_m
\ge\frac1{Q_{2^n}}\sum_{m=2^{n-1}+1}^{2^n} Q_m
\ge 2^{n-1}\,\frac{Q_{2^{n-1}}}{Q_{2^n}}\ge\frac14\,2^n.
\]

With the uniform lower bound $F_{2^n}(x\dotplus t) \ge c\,2^n$ established for 
sufficiently large~$n$ in both cases, the necessity of the $A_p(\mathcal{D})$ condition 
follows standardly.  
For $1<p<\infty$, testing hypothesis~(ii) for $t_{2^n}$ 
on $f=\chi_{I_{n,k}}w^{-1/(p-1)}$, via a standard truncation argument to ensure 
$L^p_w$-integrability, yields $w\in A_p(\mathcal{D})$. 
The case $p=1$ follows as in Theorem~\ref{th:Fejer-A1}.
\end{remark}

Consequently, as in Section~\ref{sec:var-Lebesgue}, the Rubio de Francia extrapolation 
from Theorem~\ref{th:RdF} allows us to extend Theorem~\ref{th:Norlund-Ap} to variable 
Lebesgue spaces.
\begin{theorem}\label{th:Norlund-var}
Let $p(\cdot)$ be an exponent with $1\le p_-\le p_+<\infty$. Suppose that $t_n$ is 
a regular N\"orlund mean generated by a sequence $\{q_k\}_{k\in\mathbb{N}}$ that is 
either non-decreasing and satisfies~\eqref{cond:non-decr}, or non-increasing and 
satisfies~\eqref{cond:non-incr}. Then the following are equivalent:
\begin{enumerate}
\item[\textnormal{(i)}] $p(\cdot)\in\mathcal{A}(\mathcal{D})$.

\item[\textnormal{(ii)}] There exists a constant $C>0$, depending only 
on~$p(\cdot)$ and $[p]_\mathcal{A}$, 
such that for all $n\in\mathbb{P}$ and 
$f\in L^{p(\cdot)}[0,1)$,
\[
\|t_{n}f\|_{p(\cdot)}\le C\,\|f\|_{p(\cdot)}.
\]
	
\item[\textnormal{(iii)}] For all $f\in L^{p(\cdot)}[0,1)$,
\[
\lim_{n\to\infty}\|t_{n}f-f\|_{p(\cdot)}=0.
\]
\end{enumerate}
\end{theorem}

We conclude by establishing the quantitative rate of approximation by N\"orlund means 
in variable Lebesgue spaces. In the classical $L^p[0,1)$ setting, this was first 
investigated by M\'oricz and Siddiqi~\cite{MS92} for special N\"orlund means. 
Recently, Areshidze and Tephnadze~\cite[Theorems~1 and~3]{AT24} established the 
following general estimates for means generated by monotonic sequences, bypassing the 
restrictive conditions of~\cite[Theorem~1]{MS92}:
\begin{theorem}
Let $f\in L^p[0,1)$ with $1\le p<\infty$. Suppose that $t_n$ is a regular N\"orlund 
mean generated by a sequence $\{q_k\}_{k\in\mathbb{N}}$, and $|n|$ is the unique
integer satisfying $2^{|n|}\le n<2^{|n|+1}$. The following estimates hold:
\begin{enumerate}
\item[\textnormal{(a)}] If $\{q_k\}_{k\in\mathbb{N}}$ is non-decreasing, then 
\[
\|t_nf-f\|_p\le \frac{18}{Q_n} \sum_{k=0}^{|n|-1}2^k q_{n-2^k}\,
\omega_p(f,2^{-k})
+12\,\omega_p(f,2^{-|n|}).
\]

\item[\textnormal{(b)}] If $\{q_k\}_{k\in\mathbb{N}}$ is non-increasing and satisfies
condition~\eqref{cond:non-incr}, then
\[
\|t_nf-f\|_p \le C\sum_{k=0}^{|n|} 2^{k-|n|}\,\omega_p(f,2^{-k}),
\]
where $C>0$ is a constant depending only on~$p$.
\end{enumerate}
\end{theorem}

The proof of the following theorem follows line-by-line the arguments of Areshidze and 
Tephnadze~\cite{AT24}, replacing the classical inequality~\eqref{ineq:Fej-Lp} with the 
variable Lebesgue space estimate established in Theorem~\ref{th:Fej-var-Lp}.
\begin{theorem}\label{th:Norlund-est}
Let an exponent $p(\cdot)\in\mathcal{A}(\mathcal{D})$ satisfy $1\le p_-\le p_+<\infty$,
and let $f\in L^{p(\cdot)}[0,1)$. Suppose that $t_n$ is a regular N\"orlund mean
generated by a sequence $\{q_k\}_{k\in\mathbb{N}}$, and $|n|$ is the unique integer 
satisfying $2^{|n|}\le n<2^{|n|+1}$. The following estimates hold:
\begin{enumerate}
\item[\textnormal{(a)}] If $\{q_k\}_{k\in\mathbb{N}}$ is non-decreasing, then 
\[
\|t_nf-f\|_{p(\cdot)}\le \frac{C}{Q_n} 
\sum_{k=0}^{|n|-1}2^k q_{n-2^k}\,\omega_{p(\cdot)}(f,2^{-k})
+C\,\omega_{p(\cdot)}(f,2^{-|n|}).
\]

\item[\textnormal{(b)}] If $\{q_k\}_{k\in\mathbb{N}}$ is non-increasing and satisfies
condition~\eqref{cond:non-incr}, then
\[
\|t_nf-f\|_{p(\cdot)}\le 
C\sum_{k=0}^{|n|}2^{k-|n|}\,\omega_{p(\cdot)}(f,2^{-k}).
\]
\end{enumerate}
In both cases, $C>0$ is a constant depending only on~$p(\cdot)$ and $[p]_\mathcal{A}$.
\end{theorem}
\begin{corollary}
Let an exponent $p(\cdot)\in\mathcal{A}(\mathcal{D})$ be such that 
$1\le p_-\le p_+<\infty$, and let $\{q_k\}_{k\in\mathbb{N}}$ be a sequence generating 
N\"orlund means $t_n$. Suppose that $\{q_k\}_{k\in\mathbb{N}}$ is either non-decreasing 
and satisfies~\eqref{cond:non-decr}, or is non-increasing and 
satisfies~\eqref{cond:non-incr}.
If $f\in\textnormal{Lip}(\alpha,p(\cdot))$ for some $\alpha>0$, then
\[
\|t_nf-f\|_{p(\cdot)}=
\begin{cases}
O(n^{-1}) & \textnormal{if }\alpha>1, \\
O(n^{-1}\log n) & \textnormal{if }\alpha=1, \\
O(n^{-\alpha}) & \textnormal{if }0<\alpha<1
\end{cases}
\quad\textnormal{as }n\to\infty.
\]
\end{corollary}

\bibliographystyle{abbrv}
\bibliography{KST-bib}

\end{document}